\documentclass[runningheads, orivec]{llncs}

\usepackage{amsmath}  
\usepackage{amssymb}  
\usepackage{mathtools} 
\usepackage{multirow}
\usepackage[flushleft]{threeparttable}
\usepackage{algorithm} 
\usepackage{algpseudocode}
\usepackage{hyperref}
\usepackage{changes}

\DeclareMathOperator{\PT}{PT}

\newcommand{\R}{\mathbb{R}}
\newcommand{\M}{\mathcal{M}}
\newcommand{\E}{\mathcal{E}}
\newcommand{\V}{\mathcal{V}}

\newcommand{\ein}{e_\mathrm{in}}
\newcommand{\eout}{e_\mathrm{out}}

\usepackage[T1]{fontenc}
\usepackage{graphicx}

\usepackage{overpic}

\begin{document}

\title{$p$-Laplacians for Manifold-valued Hypergraphs}
%
%

\author{Jo Andersson Stokke\inst{2} \and
Ronny Bergmann\inst{2} \and
Martin Hanik\inst{1,3,4} \and
Christoph von Tycowicz\inst{4}}

\authorrunning{J. Stokke et al.}
%

\institute{BIFOLD---Berlin Institute for the Foundations of Learning and Data, Ernst-Reuter Platz 7, 10587 Berlin, Germany \and Norwegian University of Science and Technology, NO-7491 Trondheim, Norway \email{jo.stokke.andersson@gmail.com, ronny.bergmann@ntnu.no} \and
Machine Learning Group, Technical University Berlin, Stra{\ss}e des 17. Juni 135, 10623 Berlin, Germany \and 
Zuse Institute Berlin, Takustra{\ss}e 7, 14195 Berlin, Germany
\email{\{hanik,vontycowicz\}@zib.de}}
\maketitle              
\begin{abstract}
Hypergraphs extend traditional graphs by enabling the representation of $N$-ary relationships through higher-order edges. Akin to a common approach of deriving graph Laplacians, we define function spaces and corresponding symmetric products on the nodes and edges to derive hypergraph Laplacians. While this has been done before for Euclidean features, this work generalizes previous hypergraph Laplacian approaches to accommodate manifold-valued hypergraphs for many commonly encountered manifolds. 

\end{abstract}

\section{Introduction}
    Today, graphs are used extensively to model phenomena across a varied set of domains, e.g.\ 
    image processing~\cite{lezoray2008nonlocal},
    social network theory~\cite{Huang2005InfluenceOL}, and molecular dynamics~\cite{pietrucci2011graph}. 
    However, graph modeling is confined to pairwise connections of nodes, which raises the concern that its use may be limited in the presence of complex multi-node relationships. Hypergraphs address this limitation by allowing for arbitrarily many nodes in their edges, possibly making them a more natural candidate for such phenomena. The hypergraph approach has led to improvements, e.g.\ in image processing~\cite{bretto2005hypergraph} and social network modeling~\cite{antelmi2021social}. 
    Further works have generalized the Laplacian to make the highly successful spectral graph theory available for hypergraphs~\cite{li2018submodular,saito2018hypergraph}. 
    Examples include here the application to the hypergraph minimal cut problem~\cite{saito2018hypergraph} as a regularization term in the classification of hypergraph nodes~\cite{jost2019hypergraph} and to derive diffusion equation flows~\cite{fazeny2023hypergraph}. Similarly to graphs, there has also been interest in the continuous family of $p$-Laplacians for hypergraphs~\cite{li2018submodular,fazeny2023hypergraph}.
    
    A type of data that is increasingly encountered in graph applications today is manifold-valued data; e.g.\ SPD matrix-valued graphs obtained by diffusion tensor imaging are used to analyze the function of the human brain~\cite{pennec2019riemannian,Bergmann_2018}, and surface graphs with sphere-valued normals were used for a shape classification of human hippocampi in~\cite{hanik2024manifold}.
    Furthermore, embedding (abstract) graphs in spaces of negative curvature---prominently, hyperbolic space---is well suited for learning hierarchical structures and classification~\cite{hanik2024manifold,peng2021hyperbolic}.
    To the best of our knowledge, no comprehensive framework for hypergraphs with features in general manifolds exists beyond a limited application in hyperbolic space~\cite{li2022hyperbolic}. Therefore,  translating non-Euclidean frameworks from graphs to hypergraphs is mostly unexplored.
    
    We present the first $p$-Laplacians for manifold-valued hypergraphs to close this gap. Thereby, this work opens the doors to novel application areas for hypergraphs that have not been considered before. 
    Using the approach of Bergmann and Tenbrinck~\cite{Bergmann_2018}, we define two different $p$-Laplacians. They are based on two generalizations of the mean hypergraph gradient by Fazeny, Tenbrinck, and Burger~\cite{fazeny2023hypergraph} and come with distinct features, the application deciding which one fits better.   
    Both $p$-Laplacians coincide with their counterparts for standard graphs and, in many interesting cases (e.g.\ for $p=2$), with their Euclidean-featured hypergraphs; furthermore, Dirichlet's principle applies with natural (semi) inner products for Euclidean features under the same conditions.
    
    As an application of the novel $p$-Laplacians, we define the hypergraph-diffusion equation. We solve it numerically for hypergraph domains and show different time evolutions and equilibria on the 2-sphere. The experiment shows qualitatively different behavior between the two hypergraph Laplacians and the graph Laplacian. Moreover, the two proposed approaches yield hypergraph Laplacians that exhibit different diffusion dynamics.\footnote{Proofs of our mathematical results can be found in the master's thesis~\cite{stokke2024hypergraph}; accompanying code is available at \url{https://doi.org/10.5281/zenodo.15102286}.}

\section{Preliminaries}
    We first briefly introduce the necessary theoretical tools; see, e.g.\ \cite{jost2019hypergraph,fazeny2023hypergraph} for more in-depth information on oriented hypergraphs and \cite{Lee2018} on Riemannian geometry.
    \subsection{Oriented Hypergraphs}
    The central element in our paper is the hypergraph, which we define now. In the definition, $\mathcal{P}(V)$ denotes the power set of $V$.
    \begin{definition}
        A weighted oriented hypergraph is a triple $OH = (\mathcal{V},\mathcal{E}, w)$, where 
        \begin{enumerate}
            \item $\mathcal{V} = \{1, \ldots, N\}$, $N \in \mathbb{N}$, represents the \textit{nodes} or \textit{vertices},
            \item $\mathcal{E} \subseteq \{ (e_{\mathrm{in}}, e_{\mathrm{out}}) \in \mathcal{P}(\mathcal{V}) \times \mathcal{P}(\mathcal{V}) \}$ are the oriented hyperedges,
            \item $w: \mathcal{E} \to \mathbb{R}_{\geq0}$ is the hyperedge weight function.
        \end{enumerate}
    \end{definition}
    We always assume $e_\mathrm{in} \cap e_\mathrm{out} = \emptyset$.
    Given a weighted oriented hypergraph $OH = (\mathcal{V},\mathcal{E}, w)$ we say that $\widetilde{OH} = (\mathcal{V},\widetilde{\mathcal{E}}, w)$ is the \textit{oppositely oriented hypergraph} if $$e = (e_\mathrm{in}, e_\mathrm{out}) \in \mathcal{E} \iff \Tilde{e} := (e_\mathrm{out}, e_\mathrm{in}) \in \widetilde{\mathcal{E}}.$$ Here, $\Tilde{e}$ is called \textit{opposite to} $e$.
    Finally, we say that $OH$ is \textit{symmetric} if $OH = \widetilde{OH}$ and $w(\Tilde{e}) = w(e)$ for all $e \in \E$.
    
    Local neighborhoods in hypergraphs will be essential to us. The hyperedge in-neighborhood set of a vertex $u \in \V$ is given by
    $$
        \mathcal{N}^{\mathrm{in}}(u) := \{(e_\mathrm{in}, e_\mathrm{out}) \in \E | u\in e_{\mathrm{in}} \}.
    $$
    
    \subsection{Differential Geometry}
    
    Let $(\M, \langle \cdot, \cdot \rangle)$ be a connected, complete Riemannian manifold with Riemannian metric $\langle \cdot, \cdot \rangle$ and induced distance $d_\M: \M \times \M \to \mathbb{R}$. We denote the Riemannian exponential and logarithm at $x \in \M$ by $\exp_x$ and $\log_x$, respectively. Since $\M$ is complete, a length-minimizing geodesic segment $\gamma_{xy}$ connects the two points for every pair $x,y \in \M$. This segment is unique if $y$ is inside the injectivity domain of $\exp_x$, and, in that case, we denote the parallel transport of a vector $X \in T_x\M$ to $T_y\M$ along $\gamma_{xy}$ by $\operatorname{PT}_{x \rightarrow y}X$.
    
    Recall that a Fr\'echet mean~\cite{pennec2019riemannian} of a set of points $x_1,\dots,x_n \in \M$ is a minimizer of the sum of squared distances $m:= \operatorname{argmin}_{x\in \mathcal{M}}\sum _{i=1}^n d_\M^{2}(x ,x_i)$; it is unique if the points are sufficiently localized.
    
\section{The Manifold Hypergraph Laplacian}

    We now propose two frameworks with Laplacians for manifold-valued functions on hypergraphs. They are based on two different definitions of hypergraph gradients.
    The two approaches differ in whether they average function values on the manifold using the Fr\'echet mean or the Euclidean average in a tangent space. For reasons that will become apparent, we refer to the former as the \textit{Fr\'echet} and the latter as the \textit{pairwise approach}. 
    Before we can define these notions, we must introduce function spaces on hypergraphs.
    
    \subsection{Function spaces on hypergraphs}
    
    Let $OH = (\V, \E, w)$ be a weighted hypergraph and $(\mathcal{M}, \langle \cdot, \cdot \rangle )$ a connected, complete Riemannian manifold. 
    We say that
    $f \colon \V \to \mathcal{M}$ is an \textit{$\mathcal{M}$-valued vertex function} and denote the set of all such functions by $\mathcal{H}\left(\V; \mathcal{M}\right)$. 
    In this paper, we assume that we work with vertex functions $f \in  \mathcal{H}(\V, \M)$ such that, for all $(e_\mathrm{in}, e_\mathrm{out}) \in \E$, the Fr\'echet mean $x_{\mathrm{in}}:= \operatorname{argmin}_{x\in \mathcal{M}}\sum _{u \in e_{\mathrm{in}}} d_\M^{2}\left(x,f(u)\right)$ of the edge input set and the (analogously given) Fr\'echet mean $x_{\mathrm{out}}$ of the output set are unique. Moreover, we assume that $\log_{f(v)}f(w)$
    exists for all $v \in e_{\mathrm{in}}$ and $w \in e_{\mathrm{out}}$.
    
    We need two slightly different domains for edge functions, one for each approach. 
    First, we define the \textit{discrete tangent bundle over the Fr\'echet means of the in-sets of $f$} as the disjoint union
    $T_f^F\mathcal{M} \colon= \dot{\cup}_{e \in \E} T_{x_\mathrm{in}}\mathcal{M}.$
    Calling a mapping 
    $$H_f \colon \E \rightarrow T_f^F\mathcal{M},\quad 
    e \mapsto H_f(e) \in T_{x_{\mathrm{in}}}\mathcal{M},$$ \textit{Fr\'echet-hyperedge function of $f$}, we equip the space $\mathcal{H}(\E; T_f^F \mathcal{M})$ of these functions with the \textit{Fr\'echet inner product} 
    $$\left\langle  H_f, G_f\right\rangle_{\mathcal{H}\left(\E ; T_f^F\mathcal{M}\right)}:=\sum_{e \in \E}\left\langle H_f(e), G_f(e)\right\rangle_{x_{\mathrm{in}}}.$$
    
    Second, we define the pairwise tangent bundle. To this end, we introduce the notation $\Pi_{i \in I} A_i$ to denote the Cartesian product of sets $A_i$ indexed by $I$, $A^N$ for a natural number $N\in \mathbb N$ to denote the Cartesian product of $N$ copies of set $A$ and $|A|$ to denote the cardinality of A.
    We define the \textit{discrete pairwise tangent bundle of $f$} as $T_f^P\mathcal{M} = \dot{\cup}_{e \in \E} (\Pi_{u \in \V} T_{f(u)}^{|\V|}\mathcal{M})$. A \textit{pairwise-hyperedge function of $f$} is then a mapping \begin{align*}
        H_f \colon \E &\rightarrow T_f^P\mathcal{M}, \quad e \mapsto H_f(e) \in \Pi_{u \in \V} T_{f(u)}^{|\V|}\mathcal{M};
    \end{align*} we denote the space of such functions by $\mathcal{H}(\E; T_f^P \mathcal{M})$ and and endow it with the \textit{pairwise inner product} 
    $$\langle  H_f, G_f\rangle_{\mathcal{H}(\E ; T_f^P \mathcal{M})}= \sum_{e \in \E} \sum_{u \in e_{\mathrm{in}}} \sum_{v \in e_{\mathrm{out}}}\langle  H_f(e)(u,v), G_f(e)(u,v)\rangle_{f(u)}.$$

    Note that for $\M = \R^d$ we have $T_f^F\mathbb{R}^d \cong \mathbb{R}^{d}$ and
    $T_f^P\mathbb{R}^d \hookrightarrow \mathbb{R}^{N\times N\times d}$ for $N = |\V|$,
    where $\cong $ and $\hookrightarrow$ mean ``linearly isomorphic to'' and ``linearly isomorphic to a subspace of'', respectively.
    
    \subsection{Gradients and Laplacians}
    
    We now introduce our novel framework. The goal is a definition that agrees with the standard definition in the case of Euclidean data and preserves as many of its fundamental properties as possible. However, as is often the case with generalizations to manifolds, multiple approaches are possible, each with different characteristics. For this reason, we derive two alternative generalizations, starting 
    with definitions of the gradients.
    \begin{definition}[Gradients]
    The \textit{Fr\'echet gradient} of a weighted oriented hypergraph $OH$ is given by
    \begin{equation*}
        \begin{aligned}
          &\nabla^F \colon \mathcal{H}\left(\mathcal{V} ; \mathcal{M}\right) \rightarrow \mathcal{H}\left(\mathcal{E} ; T_f^F\mathcal{M}\right) \\
          & \nabla^F f(e) :=  \sqrt{w(e)}\log_{x_{\mathrm{in}}}x_{\mathrm{out}}  
        \end{aligned}
    \end{equation*}
    and the \textit{pairwise gradient} by
    \begin{equation*}
        \begin{aligned}
              &\nabla^P \colon \mathcal{H}\left(\V ; \mathcal{M}\right) \rightarrow \mathcal{H}\left(\E ; T_f^P\mathcal{M}\right) \\
              &\nabla^P f(e)(u,v) := \frac{\sqrt{w(e)}}{|e_{\mathrm{in}}||e_{\mathrm{out}}|} \log_{f(u)}f(v).
        \end{aligned}
    \end{equation*}
    \end{definition}
    
    These definitions are motivated by the following properties, which are the analogs of fundamental characteristics of the continuous gradient.
    \begin{lemma}
        The Fr\'echet gradient has the following two properties:
        \begin{enumerate}
            \item[(i)] Locally constant vertex functions have vanishing gradients, i.e.\ if $f(u) = f(v)$ for all $u \in \ein ,v \in \eout$ in a hyperedge $e \in \E$, then $\nabla^Ff(e)=0$.
            \label{frechet grad kernal}
            \item[(ii)] If the edge weights are symmetric, then it is antisymmetric under parallel transport: 
            $$\nabla^Ff(e)= -\mathrm{PT}_{x_{\mathrm{out}} \rightarrow x_{\mathrm{in}}}\nabla^Ff (\Tilde{e})$$ for $\Tilde{e}$ opposite to $e$.
            \label{antisymgrad}
        \end{enumerate}
        The pairwise gradient has the following two properties:
        \begin{enumerate}
            \item[(iii)] It has a locally trivial kernel, i.e.\ $\nabla^P f(e) = 0$ for a hyperedge $e = (\ein,\eout) \in \E$ iff $f(u) = f(v)$ for all $u\in \ein$ and $v \in \eout$.
            \item[(iv)] If the edge weights are symmetric, then its components are pairwise antisymmetric under parallel transport: 
            $$\nabla^P f(e)(u,v) = -\PT_{f(u) \rightarrow f(v)} \nabla^P f(\Tilde e)(v,u)$$ 
            for $\Tilde{e}$ opposite to $e$.
        \end{enumerate}
    \end{lemma}\label{lemma pairwise grad prop}
    Note that $(iii)$ is stronger than $(i)$. Indeed, non-constant functions also lie in the kernel of the Fr\'echet gradient if, for all edges, $x_{\mathrm{in}} = x_{\mathrm{out}}$. Depending on the application, one of the two behaviors might be favorable.
    
    We can now use these gradients to define Laplacians for manifold-valued hypergraphs.
    
    \begin{definition}[$p$-Laplacians]
        Let $OH = (\V, \E, w)$, $p \in (0, \infty)$, and $\eta \in \{0,1\}$.
        The isotropic Fr\'echet $p$-Laplacian 
        of $OH$ is given by
        \begin{equation*}\label{isolapformula} 
            \begin{aligned}
                \Delta_p^{\mathrm{i},F} f(u)
                & := \frac{-1}{|\mathcal{N}^\mathrm{in}(u)|^{\eta}}\left(\sum_{e \in \mathcal{N}^\mathrm{in}(u)}\frac{w(e)\, d_\M^2(x_{\mathrm{in}},x_{\mathrm{out}})}{|e_{\mathrm{in}}|}\right)^{\frac{p-2}{2}} \\
                & \quad \sum_{e \in \mathcal{N}^\mathrm{in}(u)}\frac{w(e)}{|e_{\mathrm{in}}|}\mathrm{PT}_{x_{\mathrm{in}} \rightarrow f(u)}  \log _{x_{\mathrm{in}}} x_{\mathrm{out}}.
            \end{aligned}
        \end{equation*}
        and the anisotropic Fr\'echet $p$-Laplacian as 
            \begin{equation*}\label{anisolapformula} 
            \begin{aligned}
                &\Delta_p^{\mathrm{a},F} f(u) := \frac{-1}{|\mathcal{N}^\mathrm{in}(u)|^{\eta}}\sum_{e \in \mathcal{N}^\mathrm{in}(u)}\frac{\sqrt{w(e)}^p\, d_\M^{p-2}(x_{\mathrm{in}},x_{\mathrm{out}})}{|e_{\mathrm{in}}|}  
                \mathrm{PT}_{x_{\mathrm{in}} \rightarrow f(u)} \log _{x_{\mathrm{in}}} x_{\mathrm{out}}
            \end{aligned}
        \end{equation*}
        The isotropic pairwise $p$-Laplacian of $OH$ is given
        \begin{equation*}
            \begin{aligned}\label{pairwise iso p-laplace}
                \Delta_p^{\mathrm{i},P} f(u) := &\frac{-1}{|\mathcal{N}^\mathrm{in}(u)|^{\eta}} \left(\sum_{e \in \mathcal{N}^\mathrm{in}(u)} \sum_{u_1 \in e_{\mathrm{in}}} \sum_{u_2 \in e_{\mathrm{out}}} \frac{w(e)\, d_\M^2(f(u_1), f(u_2))}{|e_{\mathrm{in}}|^2|e_{\mathrm{out}}|} \right)^{\frac{p-2}{2}} \\
                & \sum_{e \in \mathcal{N}^\mathrm{in}(u)}\frac{w(e)}{|e_{\mathrm{in}}|^2|e_{\mathrm{out}}|} \sum_{u_1 \in e_{\mathrm{in}}}  \mathrm{PT}_{f(u_1) \rightarrow f(u)} \sum_{u_2 \in e_{\mathrm{out}}}\log_{f(u_1)}f(u_2).
            \end{aligned}
        \end{equation*}
        and the anisotropic pairwise $p$-Laplacian as
        \begin{equation*}\label{pairwise aniso p-laplace}
            \begin{aligned}
                \Delta_p^{\mathrm{a},P}f(u) :&= \frac{-1}{|\mathcal{N}^\mathrm{in}(u)|^{\eta}} \sum_{e \in \mathcal{N}^\mathrm{in}(u)}\frac{\sqrt{w(e)}^p}{|e_{\mathrm{in}}|^p |e_{\mathrm{out}}|^{p-1}}  \\
                & \quad  \sum_{u_1 \in e_{\mathrm{in}}}\sum_{u_2 \in e_{\mathrm{out}}}d_\M^{p-2} (f(u_1), f(u_2))\; \mathrm{PT}_{f(u_1) \rightarrow f(u)} \log_{f(u_1)}f(u_2).
            \end{aligned}
        \end{equation*}
    \end{definition}
    Setting $\eta=1$ entails normalizing the Laplacian by in-degree; this is frequently done to ensure stability and independence of sampling density~\cite{BAUER20124193}.
    The Fr\'echet and pairwise Laplacians coincide with the graph Laplacian from~\cite{Bergmann_2018} for manifold-valued graphs. Furthermore, a short calculation reveals that the isotropic and anisotropic Laplacians coincide for $p=2$. For $\eta = 0$, we thus get the 2-Laplacians
    $$
                \Delta_2^F f(u) = -\sum_{e \in \mathcal{N}^\mathrm{in}(u)} \frac{w(e)}{|e_{\mathrm{in}}|}\; \mathrm{PT}_{x_{\mathrm{in}} \rightarrow f(u)} \log_{x_{\mathrm{in}}} x_{\mathrm{out}}
    $$
    and
    $$
                \Delta_2^P f(u)= -\sum_{e \in \mathcal{N}^\mathrm{in}(u)}\frac{w(e)}{|e_{\mathrm{in}}|^{2}|e_{\mathrm{out}}|} \sum_{u_1 \in e_{\mathrm{in}}} 
                \quad \mathrm{PT}_{f(u_1) \rightarrow f(u)} \sum_{u_2 \in e_{\mathrm{out}}}\log_{f(u_1)}f(u_2).
    $$
    
    The following highlights an important property of the Fr\'echet 2-Laplacian:
    \begin{proposition}\label{laplacian is -2 standard R laplace}
        Let 
        $$\left\langle f, g\right\rangle_{\mathcal{H}\left(\V;\mathbb{R}^d\right)}^{\mathrm{stnd}} := \sum_{u \in \V} \left\langle f(u), g(u)\right\rangle_{\mathbb{R}^d}$$ 
        denote the canonical vertex inner product~\cite{saito2018hypergraph,fazeny2023hypergraph} for $\M = \R^d$. 
        Then, for all $f \in \mathcal{H}\left(\V ;\mathbb{R}^d\right)$,  
        $$\left\langle\nabla^{F} f, \nabla^{F} f\right\rangle_{\mathcal{H}\left(\E ;T_f^F\mathbb{R}^d\right)}= -\frac{1}{2}\left\langle f, \Delta_2^F f\right\rangle_{\mathcal{H}\left(\V ;\mathbb{R}^d\right)}^{\mathrm{stnd}}.$$ 
    \end{proposition}
    Dirichlet's principle thus also holds for Euclidean hypergraphs in the Fr\'echet framework: The proposition shows that the minimum of the Dirichlet energy $f \mapsto \langle \nabla^{F} f, \nabla^{F} f \rangle_{\mathcal{H}(\E; T_f^F\mathbb{R}^d)}$ is characterized by $\Delta_2^F f = 0$.
    
    With
    $$\operatorname{PTsum_u}(H_f(e)):= \sum_{u_1 \in e_{\mathrm{in}}} \sum_{u_2 \in e_{\mathrm{out}}} \mathrm{PT}_{f(u_1) \rightarrow f(u)} H_f(e)(u_1,u_2),$$ Prop.~\ref{laplacian is -2 standard R laplace} also holds for the pairwise 2-Laplacian if we replace the pairwise inner product with the semi inner product
    $$\left\langle  H_f, G_f\right\rangle_{\mathcal{H}\left(\E ; T_f^P \mathcal{M}\right)}^{\mathrm{semi}}= \sum_{u \in \V} \sum_{e \in \mathcal{N}^\mathrm{in}(u)} \frac{1}{|\ein|} \left\langle \operatorname{PTsum_u}(H_f(e)), 
    \operatorname{PTsum_u}(G_f(e))\right\rangle_{f(u)}.$$
    This can be seen by writing out $\langle\nabla^{P} f, \nabla^{P} f \rangle_{\mathcal{H}(\E; T_f^P\mathbb{R}^d)}^{\mathrm{semi}}$, and using the facts that the Euclidean parallel transport and Fr\'echet mean coincide with the identity map and the standard mean, respectively.
    
    For general manifolds, it is more nuanced: While it can be shown that $\Delta_2^P f = 0$ is a necessary minimization criterion for hypergraphs whose in-sets all have cardinality one~\cite[Sec.~4.2]{stokke2024hypergraph}, the situation is unknown for the Fr\'echet Laplacian.

\section{Experiment}
    We examine diffusion dynamics of the normalized Fr\'echet and pairwise $2$-Laplacian on unit-weight hypergraphs with features on the $2$-sphere $\M = \mathbb{S}^2$ with its canonical metric. (Computations on the sphere were performed with Morphomatics~\cite{Morphomatics}.)
    To this end, we randomly generate signed incidence matrices representing symmetric oriented hypergraphs and embed them on an octant of $\mathbb{S}^2$ by sampling angles $\theta, \phi$ from the uniform distribution on $[0, \frac{\pi}{2}) \times [0, \frac{\pi}{2}]$.
    With $\Delta_2^*\in \{\Delta_2^{F},\Delta_2^{P}\}$, we compute equilibria of the heat diffusion equation 
    $$\partial_tf(u,t) = -\Delta_2^* f(u,t), \quad u \in \V, \quad t \ge 0,$$
    using the generated embeddings as initial conditions. (See~\cite{Bergmann_2018,hanik2024manifold} for more on heat diffusion of manifold-valued graphs.) We apply the iterative numerical scheme from~\cite{Bergmann_2018} to compute the solution $f$ up to a sufficiently large $t$. 
    
    To investigate modeling differences, we also expand each generated hypergraph to a graph in the following way: Using the same (embedded) nodes, we create for each $e \in \E$ and all $u \in \ein$, $v \in \eout$ the edge $(u,v)$. We then compute an equilibrium for the graph diffusion, just like for hypergraphs. Since both of our Laplacians reduce to the one for graphs from~\cite{Bergmann_2018}, we obtain the same dynamics and equilibria. 
    
    
    \begin{figure}[ht!]
        \centering


        \includegraphics[width=\textwidth]{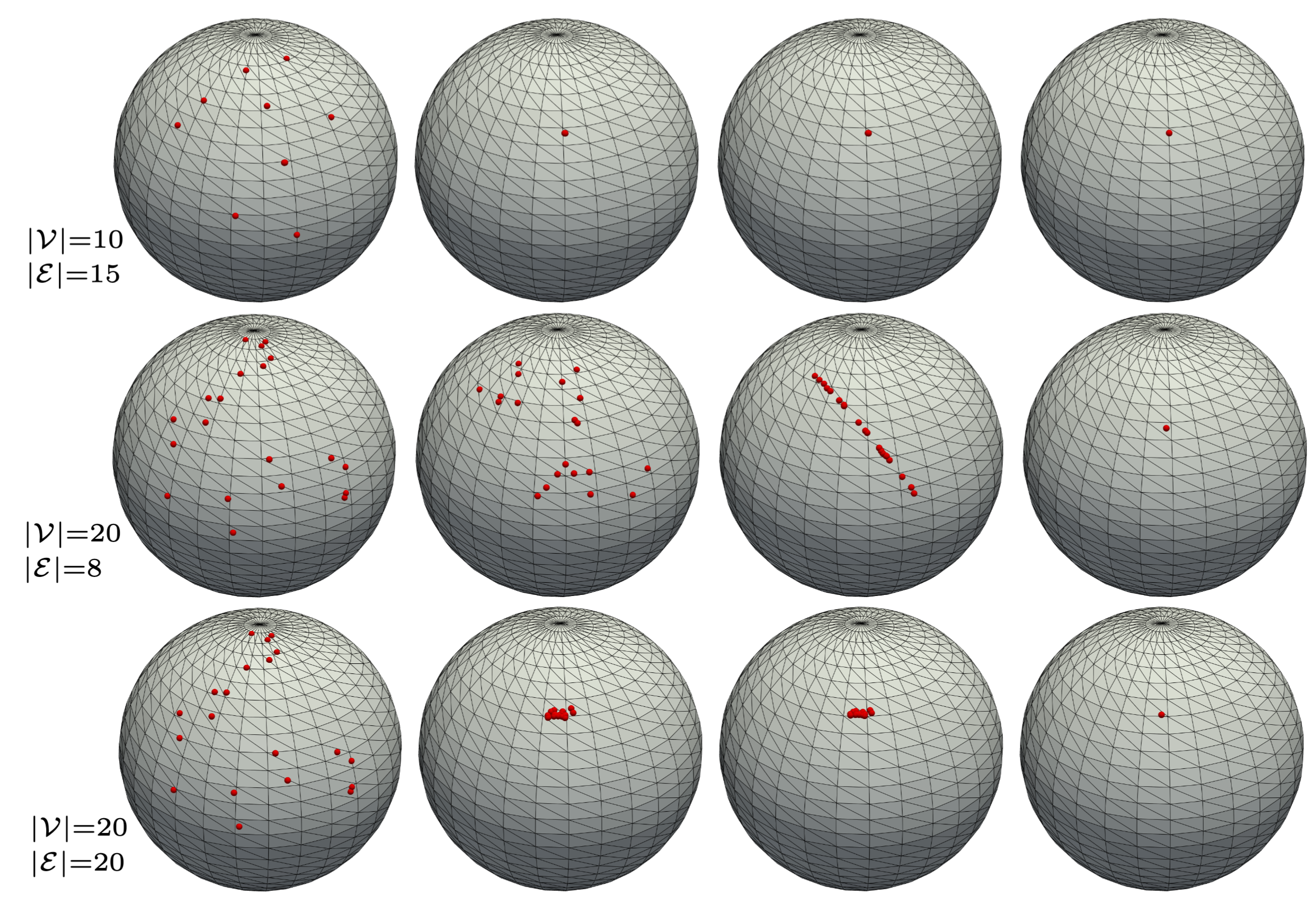}
        \caption{
        Random initial hypergraphs (first column) along equilibria of heat diffusion for $\Delta_2^{F}, \Delta_2^{P}$, and graph Laplacian (columns 2, 3, and 4, respectively).
        }
        \label{fig:results S2}
    \end{figure}
    
    The vertices of the initial graphs and the equilibrium configurations are shown in Fig.~\ref{fig:results S2}. For the first hypergraph and all graph configurations (Row 1 and Column 4), we see the typical behavior of heat diffusion processes, i.e.\ convergence to a constant value. The other two hypergraph configurations behave differently, though: Columns 2 and 3 show that the Fr\'echet mean and pairwise Laplacians exhibit equilibria other than a constant function. This interesting fact seems to occur more frequently when there are a few high cardinality edges than when there are many small ones. Additionally, even though the Fr\'echet and pairwise Laplacian tend to display similar behaviors, we see a noticeable discrepancy in the second row, where the former converges fast to a ``scattered'' solution and the latter to a configuration along a geodesic. 
    

    These results exemplify that the hypergraph Laplacians exhibit a richer set of equilibria.
    While we expect the graph Laplacian for a positively-weighted graph embedded inside a so-called \textit{normal convex neighborhood}\footnote{An open neighborhood $U \subseteq \mathcal M$ is called normal convex if all $p,q \in U$ are connected by a geodesic that does not leave $U$.}---like an octant of the sphere---to vanish only for constant configurations, we observe convergence to non-constant solutions with hypergraph.

\section{Conclusion}
    We introduced two frameworks for manifold-valued hypergraphs, each with a notion of gradient and $p$-Laplacian that generalize the case of standard graphs. The new notions conserve fundamental features of their continuous counterparts, making them canonical tools to study manifold-valued hypergraphs. Investigating the dynamics of a heat diffusion process and its equilibria, we showed that they can model a richer array of behaviors than the standard graph modeling tools.
    They thus open up a new field of applications that could not be tackled before.
    
    As there exist multiple ways of constructing a hypergraph from a graph, e.g.\ using implicit constructions like $k$-nearest neighbors \cite{dai2023hypergraph}, a natural avenue for future work is to apply hypergraph Laplacians instead of graph Laplacians when hypergraph-modeling seems promising. 
    The introduced framework is therefore directly applicable to many scenarios.
    As a first example, one can consider the various applications in \cite{Bergmann_2018} and generate different hypergraph computer vision problems (the first of which is explored in~\cite[Sec.~5.2]{stokke2024hypergraph}).  Another promising avenue is to extend the graph neural networks from \cite{hanik2024manifold} by deriving novel layers for manifold-valued hypergraphs (see~\cite[pp.~75--76]{stokke2024hypergraph}). 

\subsubsection{Acknowledgments.}  
    Funded by the Deutsche Forschungsgemeinschaft (DFG, German Research Foundation) under Germany´s Excellence Strategy---The Berlin Mathematics Research Center MATH+ (EXC-2046/1, EXC-2046/2, project ID: 390685689). 

\bibliographystyle{splncs04}
\bibliography{ref.bib}

\end{document}